\documentclass{amsart}
\usepackage{amsmath,amsfonts,amssymb,amscd,bm}
\usepackage[colorlinks=true]{hyperref}
\usepackage{tikz}
\usepackage{pgflibraryarrows}
\usepackage{pgflibrarysnakes}
\hypersetup{urlcolor=blue, citecolor=blue}
\numberwithin{equation}{section} 

\newtheorem{theorem}{Theorem}[section]
\newtheorem{corollary}[theorem]{Corollary}

\newtheorem{lemma}[theorem]{Lemma}

\theoremstyle{definition}
\newtheorem{definition}{Definition}[section]
\newtheorem{remark}{Remark}[section]

\newcommand{\beeq}{\begin{equation}}\newcommand{\eneq}{\end{equation}}
\newcommand{\al}{\alpha}    \newcommand{\be}{\beta}
\newcommand{\de}{\delta}

\newcommand{\R}{\mathbb{R}}
\newcommand{\N}{\mathbb{N}}

\newcommand{\ep}{\varepsilon}
\newenvironment{prf}{\noindent {\bf Proof.} }{\endprf\par}
\def \endprf{\hfill  {\vrule height6pt width6pt depth0pt}\medskip}
\newcommand{\pt}{\partial_t}\newcommand{\pa}{\partial}
\newcommand{\les}{{\lesssim}}
\newcommand{\supp}{\,\mathop{\!\mathrm{supp}}}

\numberwithin{equation}{section}
\title
      {Concerning ill-posedness for semilinear wave equations}
\author{Mengyun Liu}
\address{School of Mathematical Sciences\\                Zhejiang University\\                Hangzhou 310027, P. R. China}
\email{mengyunliu@zju.edu.cn}

\author{Chengbo Wang}\address{School of Mathematical Sciences\\                Zhejiang University\\                Hangzhou 310027, P. R. China}\email{wangcbo@zju.edu.cn}
\urladdr{http://www.math.zju.edu.cn/wang}

\date{\today}
\dedicatory{} \commby{}
\begin{document}
\begin{abstract}
In this paper, we investigate the problem of optimal regularity for derivative semilinear wave equations to be locally well-posed in $H^{s}$ with spatial dimension $n \leq 5$. We show this equation, with power $2\le p\le 1+4/(n-1)$, is (strongly) ill-posed in $H^{s}$ with $s = (n+5)/4$ in general. Moreover, when the nonlinearity is quadratic we establish a characterization of the structure of nonlinear terms in terms of the regularity. As a byproduct, we give an alternative proof of the failure of the local in time endpoint scale-invariant $L_{t}^{4/(n-1)}L_{x}^{\infty}$ Strichartz estimates. Finally, as an application, we also prove ill-posed results for some semilinear half wave equations.
\end{abstract}



\maketitle

\section{Introduction}
In this paper, we are interested in the problem of optimal regularity for derivative semilinear wave equations to be locally well-posed in $H^{s}$. More precisely, we consider the following Cauchy problem
\begin{equation}
\label{1.0}
\begin{cases}
\Box u = \sum_{|\gamma|=p}C_{\gamma}(\pa u)^{\gamma}, \ p\geq 2, \ p\in \mathbb{N} \\
u(0,x)=f \in H^{s}(\R^{n}), \ u_{t}(0,x)=g \in H^{s-1}(\R^{n})\\
\end{cases}
\end{equation}
where $\gamma$ is multi-index $\gamma = (\gamma_{0}, \cdot\cdot\cdot, \gamma_{n})$, $\pa = (\pa_{t}, \pa_{x_{1}}, \cdot\cdot\cdot, \pa_{x_{n}})$, $\Box= \pa_{t}^{2}-\Delta$.
A natural question is what is the minimal $s$ for which (\ref{1.0}) is locally well-posed in $H^{s}$.

Before proceeding, we recall the definition of local well-poedness.
\begin{definition}
\label{defi}
We say that the Cauchy problem (\ref{1.0}) is locally well-posed (abbreviated LWP) in $H^{s}$, if\\
(WP1) Local existence and uniqueness: given $(f, g) \in H^{s}(\mathbb{R}^{n})\times H^{s-1}(\mathbb{R}^{n})$, there exists $T =T(f,g)\in (0,\infty]$ and a unique solution $u=u(f,g)\in C([0,T);H^{s}) \cap C^{1}([0,T);H^{s-1}) \cap X^{s}_{T}$ (where $X^{s}_{T}$ is a suitable Banach space).\\ 
(WP2) The solution $u(f,g)$ depends continuously on the initial data $(f,g)$ in the following sense: for any $T_{1} < T$ and $\delta >0$, there exists $\ep >0$, such that if $\|(f-\tilde{f}, g-\tilde{g})\|_{H^{s}\times H^{s-1}}\leq \ep$, then $u(\tilde{f},\tilde{g})$ exists up to $T_{1}$ and $$\|u(f,g)-u(\tilde{f},\tilde{g})\|_{C([0,T_{1}];H^{s}) \cap C^{1}([0,T_{1}];H^{s-1})}\leq \delta.$$
(WP3) Persistence of higher regularity: if the initial data have some additional Sobolev regularity $(f, g)\in H^{\sigma} \times H^{\sigma-1}$, where $\sigma \gg s$, then the solution exists on $[0,T)$ and is in the space $C([0,T); H^{\sigma}) \cap C^{1}([0,T); H^{\sigma-1})$. In particular, if the data are in $C_{0}^{\infty}$, then the solution is smooth.\par
\end{definition}
It is known that for wave equations, any solutions obtained in a LWP framework satisfy finite speed of propagation (as is typical for the classical solutions): Let $u_{j}$ be the weak solution with initial data $(f_{j},g_{j})$ in some space-time region $\Omega_{j}\subset \R_{+}\times \R^{n}$, $j=1,2$. If $(f_{1}, g_{1}) = (f_{2}, g_{2})$ in $B_{R}(x_{0})=\{x\in \R^{n}; \|x-x_{0}\|<R\}$, then $$u_{1} = u_{2} \ \mathrm{in}\ \mathcal{D'}(\Omega),$$
where $\Omega=\Omega_{1}\cap\Omega_{2}\cap\{(t,x); \|x-x_{0}\|+t< R \}$.

Typically speaking, we say the problem is ill-posed if the problem is in contrast to the meaning of LWP.
 For wave equations, we sometimes can show the ill-posedness in strong sense (we refer it as strongly ill-posed): for any $ \varepsilon >0$, we can find initial data $(f,g)$ with $\|f\|_{H^{s}}+\|g\|_{H^{s-1}} \leq \varepsilon$ such that either it does not satisfy finite speed of propagation or it does not have local solution $u \in C([0,T];H^{s})\cap C^{1}([0,T];H^{s-1})$ for any $T > 0$.\par
Turning to the problem at hand. It is easy to see this equation is scaling invariant and the corresponding critical regularity is $s_{c}=\frac{n}{2}+1-\frac{1}{p-1}$, which is well-known to be a lower bound of $s$ that (\ref{1.0}) is locally well-posed in $H^{s}$, see e.g., Fang and Wang \cite{fw}, Tao \cite[Chapter 3]{Tao}. On the other hand, heuristically, invariance under the Lorentz transform and scaling yields another regularity  index, $$s_{l}=\frac{n+5}{4},$$ which should be another lower bound of $s$, see, e.g., \cite{fw}, where they showed (\ref{1.0}) is ill-posed in $\dot{H}^{s}$ for $s\in (s_{c}, s_{l})$ when $n= 3, 4$ and $2\leq p< 1+\frac{4}{n-1}$. Here, the result for $n=3$ and $p=2$ has been known from the work of Lindblad \cite{lbd3}. In addition, when $p=2$, it is known that the problem is LWP in $H^s$ if
$$s > \max (\frac{n}{2},\frac{n+5}{4}).$$
 When $n=2,3$, it can be proved by Strichartz estimates, see Ponce and Sideris \cite{ps} for $n=3$ and Zhou \cite{zhou} for $n=2$. In dimension $n = 4$ it can be proved within the framework of the $X^{s,b}$ spaces, see Zhou \cite{zhou}. In dimension $n \geq 5$ it was showed by Tataru \cite{dt} with modifications of $X^{s,b}$ spaces. For $p\geq 3$, the LWP results can be obtained by Strichartz estimates. Overall, the problem is LWP in $H^{s}$ for
\beeq
\begin{cases}
s\geq s_{c}, p>5, n=2,\\
s\geq s_{c}, p>3, n\geq 3,\\
s>\max(s_{c}, s_{l}), 2\leq p\leq 1+\max(2,\frac{4}{n-1}),\\
\end{cases}
\eneq
 see, e.g., Fang and Wang \cite{fw}, \cite{fw2} and references therein. Furthermore, the results can be improved if the initial data have radial symmetry or certain amount of angular regularity. For example, Fang and Wang \cite{fw2} showed that (\ref{1.0}) is LWP in $H^{s}$ with $s\geq s_{c}$ and a slight angular regularity when $n=2$, $p\geq 4$ and $s>s_{c}$ when $n=2$, $p=3$. When $n=3$ and $p=2$ (\ref{1.0}) is LWP in $H^{s}$, $s>3/2=s_{c}$ with radical symmetry data which was proved in Hidano-Jiang-Lee-Wang \cite{hjlw}. See also Sterbenz \cite{ster}, Machihara-Nakamura-Nakanishi-Ozawa \cite{mnno} and references therein for some related works.
 
 Clearly, when $2\leq p \leq 1+ \frac{4}{n-1}$, there is a gap between LWP and ill-posedness: what is the situation when $s=s_{l}$? When $n=3, p=2$ Lindblad \cite{lbd} showed that the model equation
 \beeq
 \label{mo}
 \Box u = (u_{t}-u_{x_{1}})^{2}
 \eneq
is ill-posed in $H^{2}$ (see also Lindblad \cite{lbd2}, Ettinger and Linblad \cite{eb} for quasi-linear case). Then it is natural to expect (\ref{1.0}) is ill-posed in $H^{s_{l}}$ in general when $2\leq p \leq1+\frac{4}{n-1}$. In following theorem, we show it is the case.
\begin{theorem}
\label{1}
Let $1 \leq n \leq 5$, $ 2\leq p \leq 1+ \frac{4}{n-1}$ and $p\in \N$. Here, when $n=1$ it is understood that $2 \leq p <\infty$. Then the Cauchy problem (\ref{1.0}) is (strongly) ill-posed in $H^{s_{l}}$ in general\footnote{In fact, we can show $\Box u = (u_{t}-u_{x_{1}})^{p}$ is strongly ill-posed or $\Box u=u^{2}_{t}$ ($n\leq 4$) is ill-posed in $H^{s_{l}}$, see
Theorems \ref{2.2} and \ref{4.2}.}.
\end{theorem}
\begin{remark}
When $n=1$, one can obtain the LWP  results when $s> 3/2$ by classical energy methods. Our result shows that for any $p\geq 2$, we can find a nonlinear term such that (\ref{1.0}) is ill-posed in $H^{\frac{3}{2}}$, which proves the sharpness of the LWP results in general.
\end{remark}
\begin{remark}
We fill the gap between LWP and ill-posedness when $2\leq p \leq 1+ \frac{4}{n-1}$. But as we discussed before, it is still open when $p=2$, $n\geq 6$, $s=s_{c}$ and $p=3$, $n\geq 4$, $s=s_{c}$.
\end{remark}
\begin{remark}
When $n=3$ and $p=3$, (\ref{1.0}) is LWP for $H^{2}=H^{s_c}=H^{s_l}$ data with slight regularity for angular variable, which was showed in \cite{mnno}.
\end{remark}
\begin{corollary}
\label{co}
Let $2 \leq n \leq 5$. Then the local in time endpoint Strichartz estimate
\begin{equation}
\label{stri}
\|\pa u\|_{L_{t}^{\frac{4}{n-1}}([0,T]; L_{x}^{\infty})} \leq C_{T}\|\pa u(0)\|_{\dot{H}^{\frac{n+1}{4}}}
\end{equation}
does not hold for homogeneous wave equation for any $T >0$.
\end{corollary}
Since (\ref{stri}) is scaling invariant, if (\ref{stri}) fails for some $T>0$ then it fails for any $T>0$. 

For homogeneous wave equation $\Box u = 0$ and $n \geq 2$, it is well-known that
\beeq
\label{k}
\|\pa u\|_{L_{t}^{q}([0,\infty)L_{x}^{\infty}} \les \|\pa u(0)\|_{\dot{H}^{\frac{n}{2}-\frac{1}{q}}},
\eneq
where $\max\{2,\frac{4}{n-1}\} < q < \infty$. See, e.g., Klainerman and Machedon \cite{km1}, Fang and Wang \cite{fw1}.
Note that the range of $q$ is sharp since Klainerman and  Machedon \cite{km} proved $n=3,q=2$ does not hold, Fang and Wang \cite{fw1} proved $n=2,q=4$ does not hold. Recently, Guo-Li-Nakanishi-Yan \cite{guo} showed that $n\geq 4, q=2$ (\ref{k}) fails. To the best of the authors' knowledge, it was not clear the local in time version of (\ref{k}) is true or not.

Noting that Theorem \ref{1} depends on the structure of the nonlinearity. It is interesting to determine the relation between regularity and the structure of nonlinear terms.  In the following, we investigate the problem for the most important quadratic case, in which case we would like to rewrite (\ref{1.0}) as
\begin{equation}
\label{3}
\begin{cases}
\Box u = C^{\alpha \beta}\partial_{\alpha}u \partial_{\beta}u\\
u(0,x)=f \in H^{s}(\R^{n}), u_{t}(0,x)=g \in H^{s-1}(\R^{n})\\
\end{cases}
\end{equation}
Here, we have used the convention that Greek indices $\al, \be$ range from $0$ to $n$ and the Einstein summation convention.
 It is known that when $(C^{\al \be}) = c (m^{\al \be})=diag(c,-c,\cdots, -c)$ for some $c\in\R$, that is (\ref{3}) satisfies Klainerman's null condition, the problem is locally well-posed in $H^{s}$ for any $s> \frac{n}{2}$ when $n\ge 2$, see Klainerman-Machedon \cite{km} ($n=3$), Klainerman-Selberg \cite{ks}($n\geq 2$). When $n=4$ and $C^{\alpha \beta}\partial_{\alpha}u \partial_{\beta}u = \frac{1}{2}u_{t}^{2}$, Fang and Wang \cite{fw} proved that (\ref{3}) is ill-posed in $H^{s}$, $s<\frac{9}{4}$. In the following theorem, we show that when $1 \leq n \leq 4$, (\ref{3}) is ill-posed in $H^{s}$ for any $s\in (s_c, s_{l})$, provided that it does not satisfy the null condition. In particular, this gives a characterization of the structure of the nonlinearity in terms of the regularity.
\begin{theorem}
\label{4}
Let $1 \leq n \leq 4 $.
The following statements are equivalent:
\begin{enumerate}
  \item (\ref{3}) satisfies the null condition;
  \item (\ref{3}) is LWP \footnote{when $n=1$ and $1/2<s<1$, with $C^{\al\be}=cm^{\al\be}$, it is understood that the equation \eqref{3} is satisfied in the divergence form $m^{\al\be}\pa_\be (e^{-cu}\pa_\al u)=0$} in $H^s$ for any $s > \frac{n}{2}$;
  \item there exists $ s \in (\frac{n}{2}, s_{l})$ such that
the problem satisfies finite speed of propagation and WP2 holds for the trivial solution.
  \item there exists $ s \in (\frac{n}{2}, s_{l})$ such that (\ref{3}) satisfies WP1 with finite speed of propagation.\end{enumerate}
\end{theorem}
Finally, as an application, we examine the ill-posedness of semilinear half wave equation with special nonlinear term
\begin{equation}
\label{half}
\begin{cases}
i\partial_{t}u - \sqrt{-\Delta}u = i|\Re u|^2  \\
u(0,x) = u_{0} \in H^{s}(\R^{n}).
\end{cases}
\end{equation}
We say the half wave equation (\ref{half}) is locally well-posed, if we replace initial data $(f,g)\in H^{s}\times H^{s-1}$ by $f \in H^{s}$ and replace $C([0,T);H^{s}) \cap C^{1}([0,T);H^{s-1})$ by $C([0,T);H^{s})$ in Definition \ref{defi}.
It has been proved that (\ref{half}) is locally well-posed when
$$s>\max\{\frac{n-1}{2}, \frac{n+1}{4}\}.$$
See, e.g., Fujiwara, Georgiev and Ozawa \cite{ozawa} for $n=2$, Dinh \cite{din} for $n\geq 2$, Hidano and Wang \cite{H-W} for $n\geq 1$ with more general nonlinear terms. Actually, by the connection of this problem with the nonlinear wave equations established in \cite[Section 6]{H-W}, it natural to infer the sharp range of $s$ should be
$$s> \max\{\frac{n-2}{2}, \frac{n+1}{4}\}.$$
In the following theorem, we verify this sharp result at least when the spatial dimension $n\le 3$, as well as the negative part for $n=4$.
\begin{theorem}
\label{half-ill}
Let $1 \leq n \leq 4$. The semilinear half wave equation (\ref{half})
is ill-posed in $H^{s}(\R^{n})$ for any $s \leq \frac{n+1}{4}$. Specifically, it can not satisfy WP2 for the trivial solution and WP3 at the same time.
\end{theorem}
This paper is organized as follows. In Section 2, we list some inequalities we shall use later. Section 3 is devoted to the proof of Theorem \ref{1} and
Corollary
\ref{co}. We constructed some initial data with the desired regularity and singularity, by exploiting extension theorems and adapting the three-dimensional functions appeared in Lindblad \cite{lbd} to the current setting.
These data are then used to show that there is a large class of 
equation (\ref{1.0}) (see Theorem \ref{8.5}) which is ill-posed in $H^{s_{l}}$, in the sense that WP2 for the trivial solution and WP3 can not hold at the same time with $s=s_l$. Moreover, we can prove
the model equation $\Box u = (u_{t}-u_{x_{1}})^{p}$ is ill-posed (in the sense that it does not satisfy WP2 with 
finite speed of propagation) as well as strongly ill-posed in $H^{s_l}$, see Theorem \ref{2.2}. In Section 4, we give the proof of Theorem \ref{4}.
We first give a proof of the LWP result under the null conditions, based on Nirenberg's example  \cite[page 45]{Kl80}. For the ill-posed part for the problem violating the null conditions, we reduced the problem to 
$\Box v = C^{00}v_{t}^{2} + C^{11}v_{x_{1}}^{2}+ 2C^{01}v_{t}v_{x_{1}}$ with 
$|C^{00} + C^{11}| + |C^{01}| > 0$, which is further reduced to the ODE $w''=(w')^2$,
by exploiting the invariance of the wave operator under the Lorentz transforms.
Moreover, the ill-posed result could be upgraded to strongly ill-posed result.
At last, we give a proof of Theorem \ref{half-ill} in Section 5. The idea is to take advantage of the close connection between equation (\ref{half}) and semilinear wave equation $\Box u = u_{t}^{2}$.
In the process, for the critical case $s=s_{l}$, we find that (\ref{3}) is ill-posed  in $H^{s_{l}}$ for a large class of matrices $(C^{\al \be})$, see Theorem \ref{4.2} and discussion below.

\section{Preliminary}
In this section we collect some inequalities we shall use later. Moreover, we use $A \les B$ to stand for $A \leq CB$ where the constant $C$ may change from line to line.
\begin{lemma}(Schauder estimate)
\label{98}
Let $V$ be a finite-dimensional normed vector space, let $f \in H^{s}_{x}(\R^{d}\rightarrow V) \cap L^{\infty}_{x}(\R^{d}\rightarrow V)$ for some $s\ge 0$. Let $k$ be the first integer greater than $s$, and let $F \in C^{k}_{loc}(V\rightarrow V)$ be such that $F(0) = 0$. Then $F(f) \in H^{s}_{x}(\R^{d}\rightarrow V)$ as well, with a bound of the form
$$\|F(f)\|_{H^{s}_{x}} \les_{F,\|f\|_{L^{\infty}}, V, s, d} \|f\|_{H^{s}}.$$
\end{lemma}

See Tao \cite[Lemma A.9]{Tao} for a proof.

\begin{lemma}
\label{2.7}
Let $n \ge 1$, $s > \frac{n}{2}$ and $-s \leq b \leq s$. Then there exists $C>0$, which depends only on $s, b, n$, such that we have
\beeq
\label{5.1}
\|fg\|_{H^{b}} \le C \|f\|_{H^{s}}\|g\|_{H^{b}}
\eneq
for any
 $f \in H^{s}$ and
 $g \in H^{b}$.
\end{lemma}
\begin{prf}
For fixed $s > \frac{n}{2}$ and $f \in H^{s}$, we treat $f$ as a linear operator.
By duality and interpolation, we need only to give the proof for the case $b=s$. The estimate with $b=s$ is classical, as $H^{s}$ is an algebra when $s>n/2$ (see, e.g., Tao \cite[Lemma A.8]{Tao}).
\end{prf}

\section{Proof of Theorem \ref{1}}
In this section, we give the proof of Theorem \ref{1} and
Corollary \ref{co}. At first, we investigate the ill-posedness for the following model equation
\begin{equation}
\label{kk}
\begin{cases}
\Box u = (u_{t}- u_{x_{1}})^{p}\ ,\\
(u(0),u_{t}(0))=(f, g) \in H^{s_{l}}\times H^{s_{l}-1}\ .
\end{cases}
\end{equation}

\subsection{About initial data}
Let 
$\al=\al(n)$ be a fixed small positive number depends on space dimension, 
we introduce the following smooth function for $x_{1}\in (0, 4)$,
\begin{equation}
\chi(x) = - \int_{0}^{x_{1}} \left(\ln \frac{6}{s}\right)^{\alpha} ds, \ x\in \R^{n}, \ x_{1}\in (0, 4)\ .
\end{equation}
More specifically, we will choose $\al$ such that
$\al\in (0, 1/2)$ for $n\le 3$ and
$\al\in
(0, 1-4/(n+1))$ for $n =4, 5$.
For convenience of presentation,
 we introduce the following notations, for $j\ge 1$,
$$B^{j}:=B_{1-2^{-j}}(e_{1})
=\{x\in \mathbb{R}^{n}; \|x-e_{1}\| < 1-2^{-j}\}
\ ,
\Lambda^{j}=\bigcup_{t\in [0, 1-2^{-j})}
B_{1-2^{-j}-t}(e_{1})\ ,
$$
$B^{n}_{1}:=B_{1}(e_{1})$,
$\Lambda=\bigcup_{t\in [0, 1)}
B_{1-t}(e_{1})$, 
where $e_{1}=(1,0,\cdots, 0)\in \R^{n}$.
In addition,
for fixed $t>0$, $\Lambda^{j}_{t}=\{x; (t,x)\in \Lambda^{j}\}$, $\Lambda_{t}=\{x; (t,x)\in \Lambda\}$.
\subsubsection{The extension of $\chi$ on $B_{1}^{n}$}
\begin{lemma}
\label{yt12}
Let $n = 1, 2$. There exists $\Psi(x) \in H^{s_{l}}(\mathbb{R}^{n})$ such that
$\Psi(x) = \chi(x)$ in $B^{n}_{1}$.
\end{lemma}
\begin{prf}
By extension Theorem 1.4.3.1 in Grisvard \cite{pg}, we aim to prove $\chi  \in H^{\frac{n+5}{4}}(B_{1}^{n})$. Observing that $\chi  \in L^{2}(B_{1}^{n})$, it suffice to show $\pa_{x_{1}}\chi
=-(\ln\frac{6}{x_{1}})^{\al}
 \in \dot{H}^{\frac{n+1}{4}}(B_{1}^{n})$. By Definition 1.3.2.1 in \cite{pg}, for $0< s <1$, $\Omega$ is an open subset of $\R^{n}$,
\begin{align*}
\|f\|_{\dot{H}^{s}(\Omega)} &=\left(\int_{\Omega}\int_{\Omega} \frac{|f(x)- f(y)|^{2}}{|x-y|^{n+2s}} dxdy\right)^{1/2}.
\end{align*}
In the following proof, we use $f(x)$ to denote $\pa_{x_{1}}\chi=-(\ln\frac{6}{x_{1}})^{\al}$. Notice that
\beeq\label{eq-chiprime}\pa_{x_{1}}f(x) = \alpha(\ln \frac{6}{x_{1}})^{\alpha-1}\frac{1}{x_{1}},\ \pa_{x_{1}}^{2}f(x) = \alpha(\ln \frac{6}{x_{1}})^{\alpha-2}\frac{1}{x_{1}^{2}}(1-\al-\ln\frac{6}{x_{1}}).\eneq

We first deal with the case $n=1$. For fixed $y$, let $x = y + z$ and exchange the order of integration, we are reduced to show
\begin{align*}
\|f\|_{\dot{H}^{\frac{1}{2}}(B^{1}_{1})}^{2} &= \int_{0}^{2}\int_{0}^{2}\frac{(f(x)-f(y))^{2}}{(x-y)^{2}}dxdy\\
&=\int_{0}^{2}\int_{-y}^{2-y}\frac{(f(y+z)-f(y))^{2}}{z^{2}}dzdy\\
&=\int_{0}^{2}\int_{0}^{2-z}\frac{(f(y+z)-f(y))^{2}}{z^{2}}dydz + \int_{-2}^{0}\int_{-z}^{2}\frac{(f(y+z)-f(y))^{2}}{z^{2}}dydz\\
&=2\int_{0}^{2}\int_{0}^{2-z}\frac{(f(y+z)-f(y))^{2}}{z^{2}}dydz < \infty.
\end{align*}

Let $g(z)=\int_{0}^{2-z}(f(y+z)-f(y))^{2}dy$ with $0< z <2$,
we claim that
\beeq
\label{eq-gz}
g(z) \les \left(\ln\frac{6}{z}\right)^{2(\al-1)}z \ .\eneq
With help of this estimate, we see that
$$\|f\|^{2}_{\dot{H}^{\frac{1}{2}}(B^{1}_{1})}
\les \int_{0}^{2}\frac{g(z)}{z^{2}}dz
 \les \int_{0}^{2}\left(\ln\frac{6}{z}\right)^{2(\al-1)}\frac{1}{z}dz 
 = \int_{\ln 3}^{\infty} t^{2(\al-1)}dt
 < \infty,$$
where we have used the fact that $2(\al-1)<- 1$, as $\al<1/2$.

The estimate  \eqref{eq-gz} for $z\ge 1/4$ is trivial,
as we have $f\in L^{2}(B^{1}_{1})$.
To prove the claim \eqref{eq-gz} for $0< z < 1/4$,
 we divide the integration into three cases: $y\in (0,  z^{2}], (z^{2},z], (z, 2-z)$, and denote the corresponding integral as $g_{1}, g_{2}, g_{3}$.

For the first case,  $ 0 < y \leq z^{2}$, we have $(f(y+z)-f(y))^{2}\leq f^{2}(y)$, as $0> f(y+z)> f(y)$. To control $g_{1}$,
we use integration by parts to get
$$
g_{1}(z)\leq\int_{0}^{z^{2}}f^{2}(y)dy
= yf^{2}(y)|_{y=0}^{z^{2}}-\int_{0}^{z^{2}}2 f f' ydy
=\left(\ln\frac{6}{z^{2}}\right)^{2\al}z^{2}+\int_{0}^{z^{2}}2\al\left(\ln\frac{6}{y}\right)^{2\al-1}dy.
$$
As $\al<1/2$,
$(\ln\frac{6}{y})^{2\al-1}$ is increasing, and
 we see that
$g_{1}(z)=\mathcal{O}((\ln\frac{6}{z^{2}})^{2\al}z^{2}+
(\ln\frac{6}{z^{2}})^{2\al-1}
z^{2})\les(\ln\frac{6}{z})^{2(\al-1)}z$.

Turning to the case $z^{2} < y < z$, by H\"{o}lder's inequality, we have
$$|f(y+z)-f(y)| = |\int_{y}^{y+z}f'(t)dt|\leq \|f'\|_{L^{\frac 3 2}(y,y+z)}z^{\frac 1 3}\les \left(\ln\frac{6}{z}\right)^{\alpha-1}\left(\frac{z}{y}\right)^{\frac 1 3},$$
and so
$$g_{2}(z)\les\  z^{\frac 2 3}\left(\ln\frac{6}{z}\right)^{2(\alpha-1)}\int_{z^{2}}^{z}\left(\frac{1}{y}\right)^{\frac 2 3}
dy
\les
\left(\ln\frac{6}{z}\right)^{2(\al-1)}z\ .
$$

For the remaining case,  $z \leq y < 2-z$, we observe from \eqref{eq-chiprime} that
$f''(t)<0$ for $t<2$ and so
$f'(t)< f'(y)$ for any $y<t<2$.
Then $$|f(y+z)-f(y)| = |\int_{y}^{z+y} f'(t)dt| \les
f'(y) z\les
\left(\ln\frac{6}{y}\right)^{\alpha-1}\frac{z}{y}\ , $$
and so $$g_{3}(z) \les \int_{z}^{2-z}\left(\ln \frac{6}{y}\right)^{2\alpha-2}\frac{z^{2}}{y^{2}}dy=
\int_{z}^{1/4}\left(\ln \frac{6}{y}\right)^{2\alpha-2}\frac{z^{2}}{y^{2}}dy+\mathcal{O}(z^{2})=w(z)+\mathcal{O}(z^{2}).$$
For $w(z)$, integration by part yields
\begin{align*}
w(z)&=\int_{z}^{1/4}\left(\ln \frac{6}{y}\right)^{2\alpha-2}\frac{z^{2}}{y^{2}}dy\\
&=-\left.\left(\ln \frac{6}{y}\right)^{2\alpha-2}\frac{z^{2}}{y}\right|_{y=z}^{1/4}+2(1-\al)\int_{z}^{1/4} \left(\ln \frac{6}{y}\right)^{2\alpha-3} \frac{z^{2}}{y^{2}}dy\\
&\leq \left(\ln \frac{6}{z}\right)^{2\alpha-2}z+\frac{2(1-\al)}{\ln24}w(z),
\end{align*}
and so we have $w(z)\les (\ln \frac{6}{z})^{2\alpha-2}z$.
Hence we obtain
$$g_{3}(z)\les w(z)+ z^{2}\les \left(\ln \frac{6}{z}\right)^{2\alpha-2}z.$$
This completes the proof of  the claim \eqref{eq-gz} and so is the proof of Lemma \ref{yt12} for $n=1$.

Turning to the case $n=2$. Similarly, let $x = y + z$ and exchange the order of integration
\begin{align*}
\|f(x)\|_{\dot{H}^{\frac{3}{4}}(B_{1}^{2})}^{2}&=
\int_{B^{2}_{1}}\int_{B^{2}_{1}}\frac{(f(x)-f(y))^{2}}{|x-y|^{\frac{7}{2}}}dxdy\\
&=\int_{B_{0}^{2}}\int_{\Omega_{z}}\frac{(f(y+z)-f(y))^{2}}{|z|^{\frac{7}{2}}}dydz,
\end{align*}
where $\Omega_{z} =  B_{1}^{2} \cap (B_{1}^{2}-z)$ depends on $z \in B_{0}^{2}=\{x\in \R^{2};|x|\leq 2\}$. Let $h(z) = \int_{\Omega_{z}}(f(y+z)-f(y))^{2}dy$, without loss of generality we consider $z_{1} >0$, then
$$h(z) \les \int_{0}^{2-z_{1}}\int_{|y_{2}|\leq \sqrt{1-(y_{1}-1)^{2}}}(f(y+z)-f(y))^{2}dy_{2}dy_{1}\les \int_{0}^{2-z_{1}}(f(y_{1}+z_{1})-f(y_{1}))^{2}y_{1}^{\frac{1}{2}}dy_{1}.$$
When $0<z_{1}<1/4$, by above, we have 
\begin{equation}
\left(f(y_{1}+z_{1})-f(y_{1})\right)^{2}\les
\begin{cases}
f^{2}(y_{1}), \ 0<y_{1}\leq z_{1}^{2}.\\
\left(\ln\frac{6}{z_{1}}\right)^{2(\al-1)}\left(\frac{z_{1}}{y_{1}}\right)^{\frac{2}{3}}, \  z_{1}^{2}< y_{1}<z_{1}.\\
\left(\ln\frac{6}{y_{1}}\right)^{2(\al-1)}\left(\frac{z_{1}}{y_{1}}\right)^{2}, \  z_{1}\leq y_{1}<2-z_{1}.\\
\end{cases}
\end{equation}
Thus by the similar argument in one dimension, we obtain when $0<z_{1}<1/4$
\begin{align*}
h(z) &\les \int_{0}^{z_{1}^{2}}\left(\ln\frac{6}{y_{1}}\right)^{2\al}y_{1}^{\frac{1}{2}}dy + \int_{z_{1}^{2}}^{z_{1}}\left(\ln\frac{6}{z_{1}}\right)^{2(\alpha-1)}\left(\frac{z_{1}}{y_{1}}\right)^{2/3}y_{1}^{\frac{1}{2}}dy
+\int_{z_{1}}^{2-z_{1}}\left(\ln\frac{6}{y_{1}}\right)^{2(\alpha-1)}\frac{z_{1}^{2}}{y_{1}^{2}}y_{1}^{\frac{1}{2}}dy\\
&\les \left(\ln\frac{6}{z_{1}}\right)^{2(\al-1)}z_{1}^{\frac{3}{2}}\les \left(\ln\frac{6}{|z|}\right)^{2(\al-1)}|z|^{\frac{3}{2}}.
\end{align*}
Hence
\begin{align*}
\|f\|^{2}_{\dot{H}^{\frac{3}{4}}(B_{1}^{2})}&\les \int_{|z|<1/4}\left(\ln\frac{6}{|z|}\right)^{2(\al-1)}\frac{1}{|z|^{2}}dz+\int_{1/4\leq |z|<2}\frac{h(z)}{|z|^{7/2}}dz\\
&\les \int_{0}^{1/4}\int_{0}^{2\pi}\left(\ln\frac{6}{r}\right)^{2(\al-1)}\frac{1}{r}d\theta dr+\mathcal{O}(1)\\
&\les \int_{\ln24}^{\infty}t^{-2(1-\al)}dt+\mathcal{O}(1)<\infty.
\end{align*}
\end{prf}
\begin{lemma}
\label{yt45}
 Let $n =3, 4, 5$. Then there exists $\Psi(x) \in H^{s_{l}}(\mathbb{R}^{n})$ such that
$\Psi(x) = \chi(x)~\mathrm{in}~ B^{n}_{1}.$
\end{lemma}
\begin{prf}
It is easy to see $\chi \in C(\bar{B_{1}^{n}})$ thus $\chi \in L^{2}(B_{1}^{n})$. 
On the one hand, for $ 1< p <\frac{n+1}{2}$, by (\ref{eq-chiprime}) we have
\begin{align*}
\|\pa^{2}_{x_{1}}\chi\|^{p}_{L^{p}(B_{1}^{n})}\les &\al^{p}\int_{0}^{2}\left(\ln\frac{6}{x_{1}}\right)^{(\al-1)p}\frac{1}{x_{1}^{p}}
(2x_{1}-x_{1}^{2})^{\frac{n-1}{2}}dx_{1}\\
\les&\int_{0}^{2}\frac{1}{x_{1}^{p-\frac{n-1}{2}}}dx_{1}<\infty,
\end{align*}
where we used the fact that $p-\frac{n-1}{2}<1$ and $(\ln\frac{6}{x_{1}})^{(\al-1)p}$ is increasing.
Moreover, if $p=\frac{n+1}{2}$, since $(1-\al)\frac{n+1}{2} > 1$ by above we have
\begin{align*}
\|\pa^{2}_{x_{1}}\chi\|^{\frac{n+1}{2}}_{L^{\frac{n+1}{2}}(B_{1}^{n})}
&\les \int_{0}^{2}\left(\ln\frac{6}{x_{1}}\right)^{(\al-1)\frac{n+1}{2}}\frac{1}{x_{1}}dx_{1}\\
&\les \int_{\ln 3}^{\infty} t^{-(1-\al)\frac{n+1}{2}}dt < \infty.
\end{align*}
Take $n =3$ and $p = 2$, then $\chi \in H^{2}(B_{1}^{3})$. Then by extension theorems in Adams- Fournier \cite{adam} page 147:
there is a linear operator $E_{2}$ mapping $H^{2}(B^{3}_{1})$ into $H^{2}(\R^{3})$, thus there exists $\Psi(x)\in
H^{2}(\R^{3})$ such that
$$\Psi(x) = \chi(x) ~~\mathrm{in} ~~ B^{3}_{1}.$$
On the other hand, when $n= 4, 5$, by (\ref{eq-chiprime}), if $1 < q <\frac{n+1}{4}$ we have
\begin{align*}
\|\pa^{3}_{x_{1}}\chi\|^{q}_{L^{q}(B_{1}^{n})}\les& \al^{q}\int_{0}^{2}\left(\left(\ln\frac{6}{x_{1}}\right)^{(\al-1)}-(1-\al)\left(\ln\frac{6}{x_{1}}\right)^{(\al-2)}\right)^{q}\frac{1}{x_{1}^{2q}}
(2x_{1}-x_{1}^{2})^{\frac{n-1}{2}}dx_{1}\\
\les &\int_{0}^{2}\left(\ln\frac{6}{x_{1}}\right)^{(\al-1)q}\frac{1}{x_{1}^{2q-\frac{n-1}{2}}}dx_{1} \  \les \int_{0}^{2}\frac{1}{x_{1}^{2q-\frac{n-1}{2}}}dx_{1} <\infty,
\end{align*}
where we have used the fact $2q-\frac{n-1}{2}<1$ and $(\ln\frac{6}{x_{1}})^{(\al-1)q}$ is increasing.
Furthermore, if $q= \frac{n+1}{4}$, since $(1-\al)\frac{n+1}{4} > 1$ by above we have
\begin{align*}
\|\pa^{3}_{x_{1}}\chi\|^{\frac{n+1}{4}}_{L^{\frac{n+1}{4}}(B_{1}^{n})}&
\les \int_{0}^{2}\left(\ln\frac{6}{x_{1}}\right)^{(\al-1)\frac{n+1}{4}}\frac{1}{x_{1}}dx_{1}\\
&\les \int_{\ln 3}^{\infty}t^{-(1-\al)\frac{n+1}{4}}dt < \infty.
\end{align*}
Then by extension theorems (\cite{adam} page 147): there exists a strong 3-extension operator $E_{n}$ for the region $B^{n}_{1}$. That is to say, the linear operator $E_{n}$ mapping $W^{k,r}(B^{n}_{1})$ into $W^{k,r}(\mathbb{R}^{n})$ for every integer $0 \leq k \leq 3$ and every $1 \leq r <\infty$. Hence, there exists $\Psi\in W^{3,q}(\R^{n}) \cap W^{2,p}(\R^{n})$ with $1 <p \leq \frac{n+1}{2}, 1 < q \leq \frac{n+1}{4}$ such that
$$\Psi(x) = \chi(x) ~~\mathrm{in} ~~ B^{n}_{1}\ .$$
By Gagliardo-Nirenberg inequality (see, e.g., Bahouri-Chemin-Danchin \cite[Theorem 2.44]{hajer})
\begin{equation}
\label{7.5}
\|\Psi\|_{\dot{W}^{\frac{n+5}{4},2}} \les \|\Psi\|^{\theta}_{\dot{W}^{2,p}}\|\Psi\|^{1-\theta}_{\dot{W}^{3,q}},
\end{equation}
where
$\theta = \frac{7-n}{4}$, $p = \frac{n+1}{2}$ and $q = \frac{n+1}{4}$ such that
\begin{equation}
\label{7.6}
\frac{n+5}{4}=2\theta+3(1-\theta),\ 
\frac{1}{2} = \frac{\theta}{p} + \frac{1-\theta}{q}\ .
\end{equation}
Thus we obtain the desired $\Psi \in H^{\frac{n+5}{4}}(\mathbb{R}^{n})$.
\end{prf}
\subsubsection{The extension of $\chi$ on $B^{j}$}
\begin{lemma}\label{thm-lemextchi}Let $1\leq n \leq 5$. For each $j\geq 1$, there exist $(f_{j}, g_{j})\in C_{0}^{\infty}(\R^{n})$ such that $f_{j}=\chi$ in $B^{j}$, $g_{j}=-\pa_{x_{1}}f_{j}$ and
$$\|f_{j}\|_{H^{s_{l}}}+\|g_{j}\|_{H^{s_{l}-1}}\le C_0 ,$$
where the constant $C_0>0$ is independent of $j$.
\end{lemma}
\begin{prf}
By Lemma \ref{yt12} and  Lemma \ref{yt45}, we see $\|\chi\|_{H^{s_{l}}(B_{1}^{n})}\leq \|\Psi\|_{H^{s_{l}}(B_{1}^{n})}< \infty$ with $1\leq n\leq 5$. Then for each $B^{j}$, it is easy to see $\chi \in C^{\infty}(\bar{B^{j}})\cap H^{k}(B^{j})\cap H^{s_{l}}(B^{j})$ for any integer $k$. Then by extension Theorem 1.4.3.1 in \cite{pg}, there exists $h_{j}\in C^{\infty}(\R^{n})\cap H^{k}(\R^{n})\cap H^{s_{l}}(\R^{n})$, 
$h_{j}(x) = \chi(x)$ in $B^{j}$,
and
$$\|h_{j}\|_{H^{s_{l}}(\R^{n})} \leq C \|\chi\|_{H^{s_{l}}(B^{j})} \leq C\|\Psi\|_{H^{s_{l}}(\R^{n})},$$
where $C$ is independent of $j$. \par
Fixes $\varphi \in C_{0}^{\infty}$ with $\supp \varphi \subset \{x\in \R^{n}; |x|\leq 4\}$ and $\varphi =1$ on $\{x\in\R^{3};|x|\leq 3\}$. Let $f_{j}= h_{j}\varphi$ and $g_{j} = -\pa_{x_{1}}f_{j}$, then $(f_{j},g_{j})\in C_{0}^{\infty}$. Moreover, we have
 $$\|f_{j}\|_{H^{s_{l}}}+\|g_{j}\|_{H^{s_{l}-1}}\les \|\varphi\|_{H^{\frac{n}{2}+1}}\|h_{j}\|_{H^{s_{l}}} \leq C_{0}\ ,$$
by Lemma \ref{2.7} with $b=n/2+1$.
 \end{prf}
 
With these functions in hand, we are ready to show ill-posedness for the model equation (\ref{kk}).
\subsection{Ill-posedness of model equation}
\begin{theorem}
\label{2.2}
Let $1 \leq n \leq 5$, $2\leq p \leq 1 + \frac{4}{n-1}$ and $p \in \N$ ($2\leq p< \infty$, when $n =1$). Then the Cauchy problem
\begin{equation}
\label{1.3}
\begin{cases}
\Box u= (u_{t}-u_{x_{1}})^{p}\\
u(0, x)=f,  u_{t}(0, x)= g \\
\end{cases}
\end{equation}
is (strongly) ill-posed in $H^{s_{l}}$. 
\end{theorem}

\subsubsection{Proof of ill-posedness}\label{sec:ilp1}
In this subsection, we prove the ill-posedness, in the sense that there is no WP2 for the trivial solution with finite speed of propagation.

Consider (\ref{1.3}) with initial data $(\ep f_{j}, \ep g_{j})\in C^{\infty}_{0}(\R^{n})$. By the classical local well-posed result in $H^{n+2}$, there exists a unique classical solution $u^{j} \in 
C([0,T_{j}]; H^{n+2})\cap C^1([0,T_{j}]; H^{n+1})\cap
C^{\infty}([0,T_{j}]\times \R^{n})$ for some $T_{j}>0$. By finite speed of propagation of classical solutions, we have $u^{j}(t,x) = u^{j}(t, x_{1})$ in $\Lambda^{j}=\{(t,x); \|x-e_{1}\|<1-2^{-j}-t, t<\min(1-2^{-j}, T_{j})\}$. Moreover, when $(t,x_{1})\in\{(t,x_{1}); 2^{-j}+t<x_{1}< 2-2^{-j}-t, t<\min(1-2^{-j}, T_{j})\}$, we obatin 
\begin{equation}
\label{1030}
(\partial_{t}-\partial_{x_{1}} )u^{j}(t,x_{1}) = \theta(t,x_{1})=\frac{2\ep(\ln \frac 6{x_{1}-t})^{\alpha}}{\big(1 -(2\ep)^{p-1}(p-1)t (\ln \frac 6{x_{1}-t}) ^{\alpha(p-1)}\big)^{1/(p-1)}},\footnote{Let $w(t,x_{1}) = (\pa_{t}-\pa_{x_{1}})u$, then (\ref{1.3}) becomes  ODE along characteristics
$$\frac{dw(t,t+x_{1})}{dt} = w^{p}, \ w(0,x_{1}) = -2\ep \chi'(x_{1})>0, x_1\in (2^{-j}, 2-2^{-j}-2t).$$
It is easy to obtain for $x_1\in (2^{-j}, 2-2^{-j}-2t)$
$$w(t, t+x_{1})=
\frac{w(0,x_1)}{(1-(p-1)tw^{p-1}(0,x_1))^{1/(p-1)}}
= \frac{2\ep |\chi'(x_{1})|}{(1-(2\ep)^{p-1}(p-1)t|\chi'(x_{1})|^{p-1})^{1/(p-1)}},$$
as long as $1-(p-1)tw^{p-1}(0,x_1) >0$.}
\end{equation}
 which will blow up when $h(t, x_{1}) = 1 - (2\ep)^{p-1}(p-1)t (\ln \frac 6{x_{1}-t}) ^{\alpha(p-1)} = 0$, that is, $$x_{1} = t+ \mu(t) = t+ 6\exp(-((2\ep)^{p-1}(p-1)t)^{-\frac{1}{\alpha(p-1)}}).$$
 We claim that $T_{j}\leq \mu^{-1}( 2^{-j})$ when $j\geq N$ for some $N>1$. In fact, since $\mu(0)=0$ and $\mu'(t)>0$, the singularity curve must intersect with $x_{1}=t+2^{-j}$ 
 at $(t,x_{1})=(\mu^{-1}(2^{-j}), \mu^{-1}(2^{-j})+2^{-j})$. Then if $\mu^{-1}(2^{-j})<\min(1-2^{-j}, T_{j})$, the singularity curve will intersect with the cone $\Lambda^{j}$.  This means the solution $u^{j}$ has singularity inside the cone $\Lambda^{j}$ which contradicted with the fact $u^{j}\in C^{\infty}(\Lambda^{j})$. Hence we must have $\mu^{-1}(2^{-j})\geq \min(1-2^{-j}, T_{j})$ for any $j\geq 1$, which means there exists $N>0$, such that $T_{j}\leq \mu^{-1}(2^{-j})$ when $j\geq N$. 
 
We claim that the maximal time of existence for the unique weak solution $w^{j}\in C([0,T];H^{s_{l}}) \cap C^{1}([0,T];H^{s_{l}-1})$
 with initial data $(\ep f_{j}, \ep g_{j})$, denoted by $T^{s_{l}}_{j}$, also satisfies $T^{s_{l}}_{j}\leq \mu^{-1}(2^{-j})$ when $j\geq N$, under the assumption of finite speed of propagation.
 
With help of this claim,  we see that $T^{s_{l}}_{j}$ goes to zero
as $j$ goes to $\infty$, since $\mu(0)=0$. But by Lemma \ref{thm-lemextchi} we have 
$$\|\ep f_{j}\|_{H^{s_{l}}}+\|\ep g_{j}\|_{H^{s_{l}-1}} \leq \ep C_{0}, \ \forall j\geq 1,$$
which shows 
the failure of WP2 for the trivial solution.

It remains to prove the claim, for which we prove by contradiction. Assume for some $j\ge N$, we have
$T^{s_{l}}_{j}> \mu^{-1}(2^{-j})$, i.e.,
 $\mu^{-1}(2^{-j})<\min(1-2^{-j}, T^{s_{l}}_{j})$.

In fact, we observe that
$$w(t,x_{1})=\int_{0}^{t}\theta(s,x_{1}+t-s)ds+\ep\chi(x_{1}+t),\ x_{1}\in(t+\mu(t), 2-t)$$
is a smooth solution to \eqref{1.3} in $\{(t,x_{1}); x_{1}\in(t+\mu(t), 2-t), t\ge 0\}$ with data
$\ep (\chi(x_{1}), -\chi'(x_{1}))$, satisfying
\beeq\label{eq-rl}(\partial_{t}-\partial_{x_{1}} )w(t,x_{1}) = \theta(t,x_{1})\ .\eneq
 Let $\Omega^{j}=\{(t,x); x_{1}>t+\mu(t), \|x-e_{1}\|<1-2^{-j}-t, t <\min(1-2^{-j},T^{s_{l}}_{j})\}$, $\Omega^{j}_{t}=\{x;(t,x)\in \Omega^{j}\}$. Then $w(t,x)\in C^{\infty}(\Omega^{j})$. As
$( f_{j}, g_{j})= (\chi(x_{1}), -\chi'(x_{1}))$ in $\Omega^{j}_{0}$,
by finite speed of propagation, we must have 
\beeq \label{11}w^{j}(t,x)=w(t,x)\ \mathrm{in}\ \mathcal{D'}(\Omega^{j}).\eneq

As $\mu^{-1}(2^{-j})<\min(1-2^{-j}, T^{s_{l}}_{j})$,  the singularity curve $x_{1}=t+\mu(t)$ will intersect with $x_{1}=t+2^{-j}$ at $(t^{*},x_{1})=(\mu^{-1}(2^{-j}), \mu^{-1}(2^{-j})+2^{-j})$. Then if we take $r=\frac{4n}{n-1}$ for $n\ge 2$ and 
 $r=2\max(p-1, 2)$
when $n=1$, by \eqref{eq-rl} and \eqref{11}, we have
\begin{align*}
&\int_{\Omega^{j}_{t^{*}}} |(w^{j}_{x_{1}} - w^{j}_{t})(t^{*})|^{r} dx
=\int_{\Omega^{j}_{t^{*}}} |(w_{x_{1}} - w_{t})(t^{*})|^{r} dx
 \\
\gtrsim& (2\ep)^{r}\int_{t^{*}+2^{-j}}^{2-2^{-j}-t^{*}} \frac{(\ln\frac{6}{x_{1}-t^{*}})^{\alpha r}}{h(t^{*}, x_{1})^{r/(p-1)}}
[(x_{1}-t^{*}-2^{-j})(2-2^{-j}-t^{*}-x_{1})]^{\frac{n-1}{2}} dx_{1}\\
\geq &\frac{(2\ep)^{r}(1-2^{-j}-t^{*})^{(n-1)/2}}{(\pa_{x_{1}}h(t^{*},t^{*}+2^{-j}))^{r/(p-1)}}(\ln\frac{6}{1-t^{*}})^{\alpha r}\int_{t^{*}+2^{-j}}^{1} (x_{1} - t^{*} -2^{-j})^{\frac{n-1}{2}-\frac{r}{p-1}}dx_{1}\\
\geq & \frac{(1-2^{-j}-t^{*})^{(n-1)/2}(\ln\frac{6}{1-t^{*}})^{\alpha r}}{((p-1)^{2}t^{*}\al(\ln 6\times2^{j})^{\al(p-1)-1}2^{j})^{r/(p-1)}}\int_{t^{*}+2^{-j}}^{1} (x_{1} - t^{*} -2^{-j})^{\frac{n-1}{2}-\frac{r}{p-1}}dx_{1}\\
=& \infty,
\end{align*}
where we have used the facts that\footnote{When $t+\mu(t)< x_{1} < 2-t$, $t<1$, we have
$$\pa_{x_{1}}h(t,x_{1})= (2\ep)^{p-1}(p-1)^{2}t\al\left(\ln\frac{6}{x_{1}-t}\right)^{\al(p-1)-1}\frac{1}{x_{1}-t} > 0, $$
$$\pa_{x_{1}}^{2}h(t, x_{1})=-\frac{(2\ep)^{p-1}(p-1)^{2}t}{(x_{1}-t)^{2}}\left(\ln\frac{6}{x_{1}-t}\right)^{\al(p-1)-2}
\left(\left(\ln\frac{6}{x_{1}-t}\right)-(1-\al(p-1))\right) \leq 0\ .$$}
$$h(t, x_{1}) = \int_{t+\mu(t)}^{x_{1}} \pa_{x_{1}}h(t,s)ds \leq (\pa_{x_{1}}h)(t,t+\mu(t))(x_{1} - t - \mu(t))\ ,$$
and
 $\frac{4n}{(n-1)(p-1)}-\frac{n-1}{2}>1$ for $2\leq p \leq 1+\frac{4}{n-1}$ with $n>1$. By Sobolev imbedding $ H^{s_{l}-1}(\R^{n})\subset L^{r}(\Omega^{j}_{t^{*}})$, $(w^{j}_{x_{1}}-w^{j}_{t})(t^{*})\not\in H^{s_{l}-1}(\R^{n})$ which is a contradiction.

\subsubsection{Failure of WP1 with finite speed of propagation}
Consider (\ref{1.3}) with initial data $(\ep\Psi, -\ep \pa_{x_{1}}\Psi)$. Assume 
there is a weak solution $u\in C([0, T_{s_{l}}); H^{s_{l}})\cap C^{1}([0, T_{s_{l}}); H^{s_{l}-1})$ for some $T_{s_{l}}>0$. Since $\Psi(x) = f_{j}(x)$ in $B^{j}$, we must have $u(t,x)=w^{j}(t,x)$ in $\Lambda^{j}$ with $t<\min(T_{s_{l}}, T_{j}^{s_{l}})$, $j\geq N$, if (\ref{1.3}) satisfies finite speed of propagation. By the same proof as in
Subsection \ref{sec:ilp1},
 we have $T_{s_{l}}\leq \mu^{-1}(2^{-j})$ for any $j\geq N$, which gives contradiction when we let $j$ goes to $\infty$.

\subsection{Ill-posedness for a large class of (\ref{1.0})}

\begin{theorem}
\label{8.5}
Let $k>0$ and $F(\pa u) = \sum_{|\gamma|=p}\tilde{C}_{\gamma}(\pa u)^{\gamma}$ be nonnegative function of $\pa u$.
Under the same condition of Theorem \ref{2.2}, the Cauchy problem
\begin{equation}
\label{8.2}
\begin{cases}
\Box u = k(u_{t}-u_{x_{1}})^{p} +F(\pa u)\\
u(0,x)=u_{0}, \pa_{t}u(0,x)=u_{1}
\end{cases}
\end{equation}
 is ill-posed in $H^{s_{l}}$. Specifically,  WP2  for the trivial solution and  WP3 can not hold at the same time with $s=s_l$.
 \end{theorem}
\begin{prf}
 Suppose we have
WP2 for the trivial solution with $s=s_{l}$, then for $\delta = 1, T=1$, there exists $\ep_{0}>0$, such that for any $(f,g)\in H^{s_l}\times H^{s_l-1}$ with
 $$\|f\|_{H^{s_{l}}}+\|g\|_{H^{s_{l}-1}} \leq \ep_{0},$$
there is a unique solution $u= u(f,g)\in C([0,1];H^{s_{l}})\cap C^{1}([0,1];H^{s_{l}-1})$ 
with $\|u\|_{C([0,1];H^{s_{l}})\cap C^{1}([0,1];H^{s_{l}-1})}\le 1$.
 
 Now, if we set $(F_j, G_j)=\ep (f_{j}, g_{j})\in C_{0}^{\infty}$
 with $\ep= \ep_0/C_0$, we see from Lemma \ref{thm-lemextchi} that
  $$\|F_j\|_{H^{s_{l}}}+\|G_j\|_{H^{s_{l}-1}} \leq \ep_{0}, \ \forall j\ge 1\ ,$$
and denote the corresponding solutions by $u^j$.
By persistence of regularity (WP3), $u^{j}\in C^{\infty}([0,1]\times \R^{n})$.
Then
as $F_{j}(x)$ depends only on $x_{1}$ in $B^{j}$ and
$F_{j} = F_{k}$ in $B^{j}$ when $j\leq k$, we
have $u^{j}(t,x)=u^{j}(t, x_{1})$ in $\Lambda^{j}$ and $u^{j} = u^{k}$ in $\Lambda^{j}$ when $j\leq k$, due to the finite speed of propagation of classical solutions.

Set $u= u^{j}$ in $\Lambda^{j}$, then $u(t,x)=u(t,x_{1}) \in C^{\infty}(\Lambda)$ satisfies
\begin{equation}
\label{6.6}
u_{tt} - u_{x_{1}x_{1}} = (\pa_{t}+\pa_{x_{1}})(\pa_{t}-\pa_{x_{1}})u\geq k(u_{t}-u_{x_{1}})^{p}, ~~~(t,x)~~\in ~~~\Lambda
\end{equation}
with data
\begin{equation}
\label{6.7}
u(0,x_{1})= \ep\chi(x_{1}), \quad \partial_{t}u(0,x_{1})= -\ep\chi'(x_{1}).
\end{equation} In particular, $u(t,x_1)\in C^\infty(\{(t,x_1); t<x_1<2-t, t\in (0,1)\})$.

Let $w(t,x_{1}) = (\pa_{t}-\pa_{x_{1}})u$, then (\ref{6.6}) becomes  ODE along characteristics
\beeq\label{11.1}\frac{dw(t,t+x_{1})}{dt} \geq k w^{p}, \ w(0,x_{1}) = -2\ep \chi'(x_{1})>0, x_1\in (0,2-2t).\eneq
It is easy to obtain for $x_1\in (0, 2-2t)$
$$w(t, t+x_{1})\geq
\frac{w(0,x_1)}{(1-k(p-1)tw^{p-1}(0,x_1))^{1/(p-1)}}
= \frac{2\ep |\chi'(x_{1})|}{(1-k(2\ep)^{p-1}(p-1)t|\chi'(x_{1})|^{p-1})^{1/(p-1)}},$$
as long as $1-k(p-1)tw^{p-1}(0,x_1) >0$.
Thus the solution of (\ref{6.6}) satisfies
\begin{equation}
\label{11.2}
(\partial_{t}-\partial_{x_{1}} )u(t,x_{1}) \geq \frac{2\ep(\ln \frac 6{x_{1}-t})^{\alpha}}{\big(1 -k(2\ep)^{p-1}(p-1)t (\ln \frac 6{x_{1}-t}) ^{\alpha(p-1)}\big)^{1/(p-1)}},
\end{equation}
which will blow up when $1 - k(2\ep)^{p-1}(p-1)t (\ln \frac 6{x_{1}-t}) ^{\alpha(p-1)} = 0$, that is, $$x_{1} = t+ \mu(t) = t+ 6\exp(-(k(2\ep)^{p-1}(p-1)t)^{-\frac{1}{\alpha(p-1)}})\ .$$
Observe that $t+\mu(t)\in (t, 2-t)$ for sufficiently small $t>0$, which gives a contradiction to $u(t,x_1)\in C^\infty(\{(t,x_1); t<x_1<2-t, t\in (0,1)\})$.\end{prf}

\subsection{Proof of Corollary \ref{co}}
If (\ref{stri}) holds for homogeneous wave equation $\Box u = 0$, then for $\Box u =F$ by Duhumel's principle together with energy estimate one has
$$\|\pa u\|_{L^\infty_T H^{\frac{n+1}{4}}} + \|\pa u\|_{L_{T}^{\frac{4}{n-1}}L_{x}^{\infty}} \les \|\pa u(0)\|_{H^{\frac{n+1}{4}}} + \|F\|_{L^{1}_{T}H^{\frac{n+1}{4}}}.$$
Here we denote $L_{T}^{q}L_{x}^{r}$ as $L_{t}^{q}([0,T];L_{x}^{r}(\R^{n}))$. 
Notice that if $F=(u_t -u_{x_1})^2$,
we have by Moser's inequality and H\"{o}lder's inequality
$$\|F\|_{L^{1}_{T}H^{\frac{n+1}{4}}} 
\les \|\pa u\|_{L_{T}^{\frac{4}{n-1}}L_{x}^{\infty}}\|\pa u\|_{L^{\infty}_{T}H^{\frac{n+1}{4}}}T^{\frac{5-n}{4}}\ .$$
Based on these two estimates
and contraction mapping argument,  it is easy to see that
the model equation \eqref{1.3}  with $p=2$ is
 locally well-posed in $H^{s_{l}}$, which is a contradiction to Theorem \ref{2.2}.
 
\section{Proof of Theorem \ref{4}}
To prove Theorem \ref{4}, we only need to show (\ref{3}) is LWP in $H^{s}$ for  any $s> n/2$ when it satisfies null condition, WP2 fails for the trivial solution for any $s\in (s_c, s_l)$ under the assumption of finite speed of propagation as well as the strongly ill-posed results in $H^{s}$ for any $s\in (s_{c}, s_{l})$, when the null condition is violated.

\subsection{Local well-posedness of (\ref{3}) under the null condition}
When (\ref{3}) satisfies null condition, 
that is, there exists $c\in\R$ such that $C^{\al \be}=cm^{\al \be}$,
the LWP results when $n\ge 2$ was known from the works
of Klainerman-Machedon \cite{km} and
 Klainerman-Selberg \cite{ks}. For completeness, we present a simple proof here, based on Nirenberg's example \cite[page 45]{Kl80}. 
Without loss of generality, we may assume $c=1$. Let $s > \frac{n}{2}$, consider the equation
\begin{equation}
\label{null}
\begin{cases}
m^{\al\be}\pa_\al(e^{-u}\pa_\be u)=0\\
u(0,x)=f \in H^{s}, u_{t}(0,x)=g\in H^{s-1}.
\end{cases}
\end{equation}
Notice that when $s>n/2$ and $s\ge 1$, the equation is equivalent to (\ref{3}) with $C^{\al \be}= m^{\al \be}$, that is
$\Box u = u_{t}^{2} - |\nabla u|^{2}$.

(WP1) Let $\tilde{f} = 1- e^{- f}, \ \tilde{g} =  e^{- f}g$, by Schauder estimate $\tilde{f}\in H^{s}$ and by Lemma \ref{2.7}, $\tilde{g}-g=(e^{-f}-1)g\in H^{s-1}$, thus $\tilde{g} \in H^{s-1}$. Since $\|f\|_{L^{\infty}} \leq C\|f\|_{H^{s}}\leq M(>0)$, then $\tilde{f} = 1 - e^{-f} \leq 1- e^{-M} < 1$. If we consider the homogeneous wave equation
\begin{align}
\label{hom}
\begin{cases}
\Box w = 0\\
w(0,x)=\tilde{f}, w_{t}(0,x)=\tilde{g}
\end{cases}
\end{align}
then there exists a unique solution $w \in C([0,\infty); H^{s})\cap C^{1}([0,\infty);H^{s-1})$. Hence for any $\varepsilon > 0$, there exists a $T_{0} \in (0,\infty)$, if $0 \leq t \leq T_{0}$, we have
 $$\|w(t) -w(0)\|_{H^{s}} \leq \ep,$$
from which, by Sobolev embedding, we get
 $$w(t) - \tilde{f} \leq \|w(t)-\tilde{f}\|_{L_{x}^{\infty}} \leq C\|w(t)-\tilde{f}\|_{H^{s}} \leq C\ep\ .$$
Then if we set $\varepsilon = \frac{1}{2C}e^{-M}$, we have
 $$w(t) \leq \tilde{f} + C\ep \leq 1-\frac{1}{2}e^{-M} < 1, \  0\leq t \leq T_{0}.$$
 Let $u = -\ln (1-w)$,  for $t \in [0, T_{0}]$. By Schauder estimate and Lemma \ref{2.7} again $ u \in C([0, T_{0}]; H^{s})\cap C^{1}([0, T_{0}];H^{s-1})$ satisfying (\ref{null}).
 
(WP2) Take $(f_{j}, g_{j}) \rightarrow (f, g)$ in $H^{s} \times H^{s-1}$. Let $\tilde{f_{j}} = 1-e^{-f_{j}}$, $\tilde{g_{j}} = e^{-f_{j}}g_{j}$, then $(\tilde{f_{j}}, \tilde{g_{j}}) \rightarrow (\tilde{f}, \tilde{g})$.
Consider homogeneous wave equation
\beeq
\Box w^{j} = 0, \ w^{j}(0) = \tilde{f_{j}}, \pa_{t}w^{j}(0) =\tilde{g_{j}},
\eneq
we have $w^{j} \rightarrow w$ in $C([0,\infty); H^{s})\cap C^{1}([0,\infty);H^{s-1})$ as $j \rightarrow \infty$. Thus for $\eta = \frac{1}{4C}e^{-M}$, there exists $N=N(\eta, T_0)>0$ such that
$$w^{j}(t)-w(t)\leq \|w^{j}(t)-w(t)\|_{L_{x}^{\infty}} \leq C\|w^{j}(t)-w(t)\|_{H^{s}}\leq C\eta,$$
for any $j\geq N$ and $t \in [0, T_{0}]$.
Then for any $j \geq N$, $t \in [0, T_{0}]$ and $x\in\R^n$, we have
$$w^{j}(t)\leq w(t)+C\eta \le 1-\frac{1}{2}e^{-M}+\frac{1}{4}e^{-M}=1-\frac{1}{4}e^{-M}<1\ .$$
 Let $u^{j} = -\ln(1-w^{j})$ with $j\geq N$, then $u^{j}$ satisfies (\ref{null}) with initial data $(f_{j}, g_{j})$ and $u^{j} \rightarrow u$ in $C([0,T_{0}];H^{s})\cap C^{1}([0,T_{0}];H^{s-1})$ when $j \rightarrow \infty$. \par
(WP3) If the initial data $(f,g) \in H^{\sigma}\times H^{\sigma-1}$ where $\sigma > s$, by the same argument in the proof of  WP1, we see the solution $u= -\ln(1-w)\in C([0,T_{0}];H^{\sigma})\cap C^{1}([0,T_{0}];H^{\sigma-1})$ since $w\in C([0,T_{0}];H^{\sigma})\cap C^{1}([0,T_{0}];H^{\sigma-1})$.

This completes the proof of 
local well-posedness in $H^s$ with $s > \frac{n}{2}$ for  (\ref{null}).

\subsection{Ill-posedness when (\ref{3}) violates the null condition}
For $Q(\partial u, \partial u) = C^{\alpha \beta}\partial_{\alpha}u\partial_{\beta}u$, since the matrix $(C^{\al\be})$ is symmetrical and Laplace operator is rotation invariant, we may take $C^{jk} = 0, j\neq k$ by rotation transform where $1\leq j, k \leq n$. Since
(\ref{3}) violates the null condition, that is, there is no $c$ such that
 $C^{\al\be} = cm^{\al\be}$, then there exists $1\leq j \leq n$ such that
$$|C^{00} + C^{jj}| + |C^{0j}| > 0,$$
and without loss of generality we may assume $j =1$ and $C^{00} + C^{11} \geq 0 \geq C^{01}$. Hence we are reduced to consider
 $$Q(\pa u) = C^{00}u_{t}^{2} + C^{11}u_{x_{1}}^{2}+ 2C^{01}u_{t}u_{x_{1}}+ \sum_{j=2}^{n}C^{jj}u_{x_{j}}^{2}+2C^{0j}u_{t}u_{x_{j}},$$
 with $C^{00} + C^{11} - 2C^{01} > 0$.\par
 Let $v(t, x_{1})= \frac{\delta^{2}}{F(\delta)}w(s)$ with $s = \frac{t-\beta x_{1}}{\delta}$, $\delta = \sqrt{1-\beta^{2}}$, $\beta \in (1/2, 1)$, $F(\delta) = C^{00} + C^{11}\beta^{2} - 2C^{01}\beta$. Since $F(0) = C^{00} + C^{11} - 2C^{01} > 0$, then there exists a $\delta_{0}<1/2$, such that $F(\delta) > \frac{1}{2}F(0) > 0$ when $\delta \in (0, \delta_{0})$. A simple calculation shows
\begin{equation}
\label{1.9}
\Box v = Q(\pa v)\eneq
is equivalent to
$w'' = (w')^{2}$.
This is a ODE and we take a class of solutions $w(s)= \ln \frac{1}{1-as}$ with $s<1/a$ with parameter $a> 0$. Then we obtain a class of special solutions to \eqref{1.9} with parameters $\delta \in (0, \delta_{0})$ and $b=\de/a>0$:
\beeq\label{eq-speslt}v(t, x_{1}) = \frac{\delta^{2}}{F(\delta)}\ln \frac{b}{b-t+\beta x_{1}}\ , t-\be x_1<b\ , \be=\sqrt{1-\de^2}\ .
\eneq
Notice that  the corresponding initial data are
$$v(0, x_{1}) = \frac{\delta^{2}}{F(\delta)}\ln \frac{b}{b+ \beta x_{1}} = f(x), \ v_{t}(0, x_{1}) = \frac{\delta^{2}}{F(\delta)}\frac{1}{b+\beta x_{1}} = g(x),\ x_1>-b/\be.$$

\subsubsection{About initial data}
\begin{lemma}
\label{1.7}
For above functions $f, g$, there exist $(\tilde{f}, \tilde{g}) \in C_{0}^{\infty}(\R^{n})$ such that
$$f = \tilde{f}, \ g = \tilde{g}~~~~\mathrm{in} ~~~~B_b,\ \supp \tilde{f}\cup \supp \tilde{g}\subset B_{2b}\ ,$$
where $B_{b}=\{x\in \R^{n}; |x|< b\}$.
Moreover, if $s>\frac{n}{2}$, there exists a constant $C>0$, independent of $b, \de\in (0,\de_0)\subset (0,1/2)$, such that  we have
\beeq
\label{initial}
\|\tilde{f}\|_{H^{s}} + \|\tilde{g}\|_{H^{s-1}} \le C
  b^{n-2s} \de^{\frac{n+5}{2}-2s}\ , 
\eneq
for any $s> n/2$ and $s\ge 1$.
\end{lemma}
\begin{prf}
It is easy to see $f\in C^{\infty}(B_{b})\cap W^{k,p}(B_{b})$ for any integer $k\ge 0$ and $1< p< \infty$. Then by extension theorem (\cite{adam}, page 147), there exists $h \in \cap_{k\ge 0, 1<p<\infty} W^{k,p}(\R^{n})\subset C^{\infty}(\R^{n})$, such that
$h=f$ in $B_{b}$.
For some fixed $\psi(x)=\phi(x/b)\in C_{0}^{\infty}$, with $\phi=1$ on $B_1$ with support in $B_{2}$, we set $\tilde{f} = h\psi\in C_{0}^{\infty}$ with $\tilde{g} =-\pa_{x_{1}}\tilde{f}/\be$.

Turning to the proof of \eqref{initial}, we first calculate the Sobolev norms for $f$ in $B_b$.
When $k \geq 1$ and $p \in (1,\infty)$ with $kp > \frac{n+1}{2}$, a simple calculation gives us
$$\|f\|_{\dot{W}^{k,p}(|x|< b)} = \frac{\de^{2}}{F(\de)}\left\|\frac{\be}{b+\be x_{1}}\right\|_{\dot{W}^{k-1,p}(|x|< b)} \les_{k,n,p}  \frac{
1}{F(0)}
b^{\frac{n}{p}-k}\de^{\frac{n+1}{p}+2-2k}\ ,
$$
where we have used the fact that 
$1-\be \sim \de^2$ and $F(\de)>\frac{1}{2}F(0)$ for any $\de\in(0, \de_{0})\subset (0, 1/2)$.
In addition,
\begin{align}
\label{nb2}
\|f\|_{L^{p}(|x|< b)} \les \|f\|_{L^{\infty}(|x|<b)}b^{\frac{n}{p}}\les \ \de^{2}b^{\frac{n}{p}}|\ln \de|.
\end{align}
Notice that there exists a constant independent of $b>0$ such that
$\|h\|_{W^{k,p}(\R^{n})}\le C \|f\|_{W^{k,p}(|x|< b)}$,
 we obtain that 
\beeq\label{nb}
\|h\|_{W^{k,p}(\R^{n})}\le C \|f\|_{W^{k,p}(|x|< b)}\le \tilde C
b^{\frac{n}{p}-k}\de^{\frac{n+1}{p}+2-2k},
\ k\geq 1, kp>(n+1)/2,
\eneq
for any $b, \de\in (0,\de_0)\subset (0,1/2)$. 

For the $H^s$ norm of $h$, we  claim that, if $s \geq\max(\frac{n}{2}, 1)$, we have
\beeq
\label{nb3}
\|h\|_{H^{s}} \les  b^{\frac{n}{2}-s}\de^{\frac{n+5}{2}-2s}.\eneq
At first, 
\eqref{nb3} is trivial for $s\ge [(n+5)/4]$ in view of \eqref{nb} with $p=2$.
Then it suffices to consider $s\in [\max(1,n/2), [(n+5)/4])$, which occur only if $n=3, 4$. However,
we know
from \eqref{nb} that
$h\in  W^{1, 2s}\cap W^{2,s}$ as long as
$s\in ((n+1)/4, 2)$.
Thus \eqref{nb3} is true for $s\in ((n+1)/4, 2)$, in view of the following complex interpolation relation
$$[W^{1, 2s}, W^{2,s}]_{s-1}=H^s\ , \ s\in (1, 2)\ .$$
This completes the proof of the claim \eqref{nb3},  by
noticing that  
 $$[\max(1,n/2), [(n+5)/4])
=[n/2, 2)
\subset((n+1)/4, 2)$$
for $n=3, 4$.

As $b<1$, we have for any $s\ge 0$,
\beeq\label{nb5}\|\psi\|_{H^{s}}
=\|\phi(\cdot/b)\|_{H^{s}}
\les b^{\frac{n}{2}-s}\ .\eneq
Then if $s>\frac{n}{2}$ and $s\ge 1$, we have
\begin{equation}
\label{nb4}
\|\tilde{f}\|_{H^{s}} +\|\tilde{g}\|_{H^{s-1}} \les
\|\tilde{f}\|_{H^{s}} =\|\psi h\|_{H^{s}} 
\les \ \|\psi\|_{H^{s}}\|h\|_{H^{s}}\les
b^{n-2s}\de^{\frac{n+5}{2}-2s}
\end{equation}
in view of
(\ref{nb3}) and (\ref{nb5}).
\end{prf}

\subsubsection{Failure of WP2 for the trivial solution for any $s\in (s_c, s_l)$}\label{sec:ilp}
Consider the Cauchy problem
\begin{equation}
\label{1.8}
\begin{cases}
\Box u = C^{00}u_{t}^{2} + C^{11}u_{x_{1}}^{2}+ 2C^{01}u_{t}u_{x_{1}}+ \sum_{j=2}^{n}C^{jj}u_{x_{j}}^{2}+2C^{0j}u_{t}u_{x_{j}}\\
u(0,x) = \tilde{f}(x), u_{t}(0,x) = \tilde{g}(x).\\
\end{cases}
\end{equation}
By the classical local well-posed result in $H^{n+2}$, there exists a unique classical solution $u \in 
C([0,T]; H^{n+2})\cap C^1([0,T]; H^{n+1})\cap
C^{\infty}([0,T]\times \R^{n})$ for some $T>0$. By finite speed of propagation of classical solutions, we obtain
$$u(t,x) = v(t, x_{1}) = \frac{\delta^{2}}{F(\delta)}\ln \frac{b}{b-t+\beta x_{1}} \quad \mathrm{in} \quad \Lambda=\{(t,x); |x|< b-t, t <\min(b,T)\}\ ,$$
which blows up at $x_1=0$ as $t\rightarrow b-$ and so $T<b$.

We claim that the maximal time of existence for the unique weak solution $w\in C([0,T];H^{s}) \cap C^{1}([0,T];H^{s-1})\subset L^\infty ([0,T]\times \R^n)$, denoted by $T_s$, also satisfies $T_s<b$, for any $s>n/2$. In fact, by finite speed of propagation, $w=u$ in $\mathcal{D'}(\Lambda)$, then we have $$\|w\|_{L^\infty([0,b)\times \R^n)}\ge
\|v(t,0)\|_{L^\infty_t([0,b))}=\infty$$
if $T\ge b $.

 For any fixed $s\in (s_c, s_l)=(n/2, ( n+5)/4)$,
 we set $s^*=\max(s, 1)\in  (s_c, s_l)$.
Let $k>\max(2\frac{s^*-s_c}{s_l-s^*}, 1)$, we choose $\de=b^k$ for any $b\in (0,\de_0)$.
 Then \eqref{initial} gives us
$$\|\tilde{f}\|_{H^{s}} + \|\tilde{g}\|_{H^{s-1}} \le \|\tilde{f}\|_{H^{s^*}} + \|\tilde{g}\|_{H^{s^*-1}}\le C  b^{n-2s^*} \de^{\frac{n+5}{2}-2s^*}
\le C b^{2(s^*-s_c)}
 \ , $$
 which tends to zero as $b$ goes to zero.
Combining it with the fact $T_s<b$, we see the failure of 
continuously dependence of the data for the trivial solution, which completes the proof.
  
\subsubsection{Proof of strong ill-posedness}
We recall
that,
in the process of the proof in Subsection \ref{sec:ilp},
we have actually constructed a series of $C_0^\infty$ data which are small in both the support and the $H^s\times H^{s-1}$ norm, while the corresponding maximal time of existence remains small. These facts could be used to boost the ill-posed result to the strongly ill-posed result.

Actually, 
for any $s\in (\frac{n}{2}, \frac{n+5}{4})$ and
 small $\ep>0$, we could construct two functions $\phi^{\star}\in H^{s}(\R^{n})\cap C^{\infty}(\R^{n}\backslash\{0\})$,
$\varphi^{\star}\in H^{s-1}(\R^{n})\cap C^{\infty}(\R^{n}\backslash\{0\})$
 with norm bounded by $\ep$  and $\supp (\phi^{\star},\varphi^{\star})\subset B_\ep (\ep e_1)$, for which there is no local solutions satisfying finite speed of propagation, with data $(\phi^{\star}, \varphi^{\star})$.

In fact, as we show above, we can find for any $j\ge 1$, a sequence $(\phi_{j}, \varphi_{j})\in C_{0}^{\infty}$,  with $\supp (\phi_{j}, \varphi_{j})\subset B_{2^{-j-1}}(3 \times 2^{-j}e_{1})$, 
$$\|\phi_{j}\|_{H^{s}} + \|\varphi_{j}\|_{H^{s-1}} \le 2^{-j},$$
and the corresponding solutions $u_{j}=u_{j}(\phi_{j},\varphi_{j})$,
as a unique weak solution in $C([0,T];H^{s}) \cap C^{1}([0,T];H^{s-1})$, exists on $[0,T_{j})$ with $T_{j} < 2^{-j-1}$.

Let $N\ge 1$ such that $2^{1-N/2}\le \ep$, 
we set $\phi^{\star} = \sum_{j=N}^{\infty}\phi_{j}$, $\varphi^{\star} = \sum_{j=N}^{\infty}\varphi_{j}$, then $(\phi^{\star}, \varphi^{\star})\in C^{\infty}(\R^{n}\backslash\{0\})$ and 
 $$\supp (\phi^{\star}, \varphi^{\star})\subset 
B_{2^{1-N/2}} (2^{1-N/2} e_1)
\subset B_\ep (\ep e_1),$$
$$\|\phi^{\star}\|_{H^{s}} +\|\varphi^{\star}\|_{H^{s-1}}\le \sum_{j=N}^{\infty} 2^{-j} \le 2^{1-N}\le \ep.$$
\par
Consider (\ref{1.8}) with initial data $(\phi^{\star}, \varphi^{\star})$, if it satisfies finite speed of propagation, then
in each ball $\Omega_{j} = \{ x: \|x-3\times2^{-j}e_{1}\| < 2^{-j-1}\}$, $(\phi^{\star}, \varphi^{\star})= (\phi_{j}, \varphi_{j})$ and we have $T_s \le T_j< 2^{-j-1}$ for any $j\ge  N$, hence $T_s = 0$, that is, there is no local solutions.

\section{Ill-posedness for semilinear half wave equation}
In this section, we consider the ill-posedness of semilinear half wave equation with special nonlinearity. The idea is to take advantage of the close connection between equation (\ref{half}) and semilinear wave equation
\begin{equation}
\Box u = u_{t}^{2}
\end{equation}
which appeared in Hidano and Wang \cite{H-W}.
\begin{theorem}
\label{4.2}
Let $1 \le n \le 4$. Then the Cauchy problem
\beeq
\label{4.1}
\begin{cases}
\Box u = u_{t}^{2}\\
u(0,x)= f, u_{t}(0,x) =g
\end{cases}
\eneq
is ill-posed in $H^{s_{l}}$, in the sense that  WP2  for the trivial solution and  WP3  can not hold at the same time with $s=s_l$.\end{theorem}
\begin{prf} 
We argue by contradiction.
Suppose we have WP2 in $H^{s_{l}}$ for the zero solution, then for $\delta =1$ and $T=1$, there exists $\ep_{1}>0$ such that if
 $$\|f\|_{H^{s_{l}}}+\|g\|_{H^{s_{l}-1}} \leq \ep_{1}$$
the solution $u\in C([0,1];H^{s_{l}})\cap C^{1}([0,1];H^{s_{l}-1})$ and
\beeq
\label{920}
\|u\|_{C([0,1];H^{s_{l}})\cap C^{1}([0,1];H^{s_{l}-1})} \leq 1.
\eneq

By Schauder estimate Lemma \ref{98}, there exists $C_1>0$ such that we have 
$$\|\ln(1-\ep f)\|_{H^{s_{l}}} 
+\|\frac{\ep \pa_{x_1}f}{1-\ep f}\|_{H^{s_{l}-1}}
\leq C_{1}(\|f\|_{L^\infty}) \ep\|f \|_{H^{s_{l}}},$$
for any $\ep\in (0, 1]$.
By Lemma \ref{thm-lemextchi},
recall that $g_{j}=-\pa_{x_{1}}f_{j}$, 
we take a sequence of initial data $(\tilde{F}_{j}, \tilde{G}_{j})=(-2\ln(1-\ep f_{j}), \frac{2\ep g_{j}}{1-\ep f_{j}})\in C_{0}^{\infty}$ with $\ep=\ep_{1}/(2C_{0}C_{1}(C_0))$, such that
 we have
$$\|-2\ln(1-\ep f_{j})\|_{H^{s_{l}}}+\|\frac{2\ep g_{j}}{1-\ep f_{j}}\|_{H^{s_{l}-1}} \leq 2\ep C_{1}(C_0)\|f_{j}\|_{H^{s_{l}}} \leq \ep_{1}\ .$$

Let $u^{j}\in C([0,1];H^{s_{l}})\cap C^{1}([0,1];H^{s_{l}-1})$ be the corresponding solution with data
$(\tilde{F}_{j}, \tilde{G}_{j})$.
By WP3, $u^j\in C^\infty([0,1]\times\R^n)$.
By the similar argument in the proof of Theorem \ref{2.2}, we know that $u^{j}(t,x)=u^{j}(t,x_{1})\in C^{\infty}(\Lambda^{j})$ and $u^{k} = u^{j}$ in $\Lambda^{j}$ when $k \geq j$.
Since $s_{l} = \frac{n+5}{4} > \frac{n}{2}$, by (\ref{920}) and Sobolev embedding there exists a constant $M > 0$ independent of $j$ such that $u^{j}(t,x) > -M$ for any $t\in [0,1]$ and $x\in \R^{n}$. \par
Let $w^{j} = 1 - e^{-\frac{u^{j}}{2}}$, then $w^{j}(t,x)= w^{j}(t,x_{1})\in C^{\infty}(\Lambda^{j})$ and $w^{j}$ satisfies
\begin{align}
\nonumber
\Box w^{j} &= \frac{1}{2}e^{-\frac{u^{j}}{2}}(\Box u^{j}+ \frac{1}{2}|\nabla u^{j}|^{2}-\frac{1}{2}(u^{j}_{t})^{2}) \nonumber
\\&= \frac{1}{4}e^{-\frac{u^{j}}{2}}( |\nabla u^{j}|^{2}+(u^{j}_{t})^{2})\label{4.3}\\ 
&\ge \frac{1}{8}e^{-\frac{u^{j}}{2}}(u^{j}_{t}-u^{j}_{x_{1}})^{2}
=\frac{1}{2}e^{\frac{u^{j}}{2}}(w^{j}_{x_{1}} - w^{j}_{t})^{2}\ge \frac{1}{2}e^{-\frac{M}{2}}(w^{j}_{x_{1}} - w^{j}_{t})^{2}\nonumber
\end{align}
inside $\Lambda^{j}$ with initial data
\beeq
\label{4.4}
w^{j}(0) = \ep\chi, \ \pa_{t}w^{j}(0) = -\ep\chi' ~~~\mathrm{in} ~~~B^{j}.
\eneq
Let $w = w^{j}$ in $\Lambda^{j}$, then $w\in C^\infty(\Lambda)$ satisfies (\ref{4.3}) in $\Lambda$ with initial data $(\ep\chi, -\ep\chi')$ in $B_{1}^{n}$. On the other hand,  by (\ref{11.1}) (\ref{11.2}) in Theorem \ref{8.5} with $k=
\frac{1}{2}e^{-\frac{M}{2}}$, we know that $w\not\in C^\infty(\Lambda)$, which gives the desired contradiction.\end{prf}

As we discussed in the introduction, the ill-posedness of equation (\ref{3}) depends on the structure of nonlinearity. And now it is clear for null condition. By the argument in Section 4, we are reduced to consider
\beeq
\label{5.2}
Q(\pa u) = C^{00}u_{t}^{2} + C^{11}u_{x_{1}}^{2}+ 2C^{01}u_{t}u_{x_{1}}+ \sum_{j=2}^{n}C^{jj}u_{x_{j}}^{2}+2C^{0j}u_{t}u_{x_{j}}.
\eneq
In fact, by the same way of Theorem \ref{4.2} the authors find that (\ref{5.2}) is ill-posed in $H^{s_{l}}$ when $|C^{00} + C^{11}| \ge 2 |C^{01}|$.
Since in this case, by applying the transform $w = 1- e^{-au}$, we can obtain
$$\Box w \geq C(w_{x_{1}} - w_{t})^{2} $$
for some constant $C > 0$, which we can handle due to Theorem \ref{8.5}.
 However, when $|C^{00} + C^{11}| < 2 |C^{01}|$, the similar argument does not work and we do not know how to handle. Anyway, based on Theorem \ref{4}, we conjecture that the problem is ill-posed in $H^{s_{l}}$ whenever the null condition is violated.
\begin{theorem}
Let $1 \leq n \leq 4$. The Cauchy problem
\begin{equation}
\label{9.6}
\begin{cases}
i\partial_{t}u - \sqrt{-\Delta}u = i|\Re u|^2 \\
u(0,x)=u_{0}\\
\end{cases}
\end{equation}
is ill-posed in $H^{s_{l}-1}$. Specifically, WP2  for the trivial solution and  WP3  can not hold at the same time with $s=s_l$.
 \end{theorem}
\begin{prf}
We suppose (\ref{9.6}) is locally well-posed in $H^{s_{l}-1}$. By  WP2, for $\delta =1, T=1$, there exists $\ep_{2}>0$, such that if
$$\|u_{0}^{j}\|_{H^{s_{l}-1}} \leq \ep_{2},$$
the solution $u^{j}=u^{j}(u_{0}^{j})$ exists up to $T=1$ and
$$\|u^{j}\|_{L^{\infty}([0,1];H^{s_{l}-1})} \leq 1.$$
 Now we take $u_{0}^{j}=\frac{2\ep g_{j}}{1-\ep f_{j}}+2 i\sqrt{-\Delta}\ln(1-\ep f_{j})$ with $\ep$ sufficiently small such that
$\|u_{0}^{j}\|_{H^{s_{l}-1}} \le  \ep_{2}$.

For each $j$, introducing $v^{j}$ by
\begin{equation}
\label{10.20}
i\partial_t v^{j} + \sqrt{-\Delta}v^{j} = i u^{j}, \quad v^{j}(0)=-2\ln(1-\ep f^{j}) .
\end{equation}
Then $\pa_{t}v^{j}(0)=\frac{1}{i}(iu^{j}(0)-\sqrt{-\Delta}v^{j}(0))=\frac{2\ep g_{j}}{1-\ep f_{j}}$ and
 $$(\pa_{t}^{2}-\Delta)v^{j}=-(i\pa_{t}-\sqrt{-\Delta})(i\pa_{t}+\sqrt{-\Delta})v^{j} =-i(i\pa_{t}u^{j}-\sqrt{-\Delta}u^{j})=(\Re u^{j})^{2}.$$
Since $f_{j}$ is real-valued, then $v^{j}$ is real-valued. So is $\sqrt{-\Delta} v^{j}$ (as $\overline{\sqrt{-\Delta}v^{j}}=
\sqrt{-\Delta}\bar v^{j}$), thus
\begin{equation}
\label{10.30}
-i\partial_t v^{j} + \sqrt{-\Delta}v^{j} = -i \bar{u^{j}}.
\end{equation}
Combine (\ref{10.20}) and (\ref{10.30}) we have
\beeq
\label{10.41}
\pa_{t}v^{j}=\frac{u^{j}+\bar{u^{j}}}{2}=\Re u^{j}, \sqrt{-\Delta}v^{j}=i\frac{u^{j}-\bar{u^{j}}}{2}=\Im u^{j}.
\eneq
Then $v^{j} \in C([0,1];H^{s_{l}})\cap C^{1}([0,1];H^{s_{l}-1})$ satisfies
\begin{equation}
\label{9.56}
\begin{cases}
(\partial_t^2-\Delta)v^{j}=(\pt v^{j})^{2}\\
v^{j}(0)=-2\ln(1-\ep f_{j}), \pa_{t}v^{j}(0)=\frac{2\ep g_{j}}{1-\ep f_{j}}.
\end{cases}
\end{equation}
In addition, there exists a constant $C>0$, independent of $\ep\in (0, 1]$, such that
$$\|v^{j}\|_{C([0,1];H^{s_{l}})\cap C^{1}([0,1];H^{s_{l}-1})}\le C,\  \forall j\geq 1\ ,$$
by which we know
$$\|v^{j}(t)\|_{L_{x}^{\infty}} \leq C_{3}, \ \forall t\in (0,1),\ j\geq 1\ ,$$
for some $C_{3}>0$.

By persistence of regularity $u^{j} \in C^{\infty}([0,1]\times \R^{n})\cap C([0,1];H^{k}(\R^{n}))$ for any $k\ge 0$, then by (\ref{10.41}), $v^{j} \in C^{\infty}([0,1]\times \R^{n})$.
 Since the initial data of the solution $v^{j}$ only depends on $x_{1}$ inside the ball $B^{j}$, then by finite speed of propagation of classical solutions, $v^{j}(t,x)=v^{j}(t,x_{1})$ in $\Lambda^{j}$. Let $w^{j} = 1- e^{-\frac{v^{j}}{2}}$, then $w^{j}\in C^{\infty}([0,1]\times \R^{n})$ and $w^{j}(t,x) =w^{j}(t, x_{1})$ in $\Lambda^{j}$. It is easy to obtain
$$\Box w^{j}=\frac{1}{4}e^{-\frac{v^{j}}{2}}\big((v^{j}_{t})^{2}+(\nabla v^{j})^{2}\big)\geq \frac{1}{2}e^{\frac{v^{j}}{2}}(w^{j}_{t}-w^{j}_{x_{1}})^{2}.$$
 Then  $w^{j}$ satisfies
\begin{align}
\label{901}
\begin{cases}
\Box w^{j}\geq \frac{1}{2}e^{-\frac{C_{3}}{2}}(w^{j}_{t}-w^{j}_{x_{1}})^{2}\\
w^{j}(0)=\ep f_{j}, w_{t}^{j}(0)=\ep g_{j}.\\
\end{cases}
\end{align}

Since $f_{j}=f_{k}$ in $B^{j}$ when $j\leq k$, by finite speed of propagation $v^{j} = v^{k}$ in $\Lambda^{j}$ when $j\leq k$, then $w^{j}=w^{k}$ in $\Lambda^{j}$ when $j \leq k$. Set $w=w^{j}$ in $\Lambda^{j}$, then $w\in C^{\infty}(\Lambda)$ satisfy (\ref{901}) inside $\Lambda$ with initial data $(\ep \chi, -\ep \chi')$ in $B_{1}^{n}$. 
However,
by (\ref{11.1}) (\ref{11.2}) in the proof of Theorem \ref{8.5} with $k=\frac{1}{2}e^{-\frac{C_{3}}{2}}$, we know that $w\not\in C^{\infty}(\Lambda)$, which is a contradiction.
\end{prf}

\subsection*{Acknowledgment}
The second author would like to thank Professor Gang Xu for helpful discussion on the extension of initial data in low dimensions.

\bibliographystyle{plain1}

\begin{thebibliography}{10}
%
\bibitem{adam}R. A. Adams and J. J. F. Fournier,
\textit{Sobolev spaces}, Second edition. Pure and Applied Mathematics (Amsterdam), 140. Elsevier/Academic Press, Amsterdam,
(2003).
%
\bibitem{hajer}H. Bahouri, J.-Y. Chemin and R. Danchin,
\textit{Fourier analysis and nonlinear partial differential equations},
Grundlehren der Mathematischen Wissenschaften [Fundamental Principles of Mathematical Sciences], 343. Springer, Heidelberg, (2011).
%
 \bibitem{din}V. D. Dinh,
 \textit{On the Cauchy problem for the nonlinear semi-relativistic equation in Sobolev spaces},
  Discrete Contin. Dyn. Syst. {\bf 38} (2018), 1127--1143.
%
\bibitem{eb} B. Ettinger and H. Lindblad,
\textit{A sharp counterexample to local existence of low regularity solutions to Einstein equations in wave coordinates},
Ann. of Math. {\bf 2} (2017), 311--330.
%
\bibitem{fw1}D. Fang and C. Wang,
 \textit{Some remarks on Strichartz estimates for homogeneous wave equation},
Nonlinear Anal. {\bf 65} (2006), 697--706.
%
\bibitem{fw}D. Fang and C. Wang,
\textit{Local well-posedness and ill-posedness on the equation of type $\Box u = u^{k}(\partial u)^{\alpha}$},
Chin. Ann. Math. {\bf 3} (2005), 361--378.
%
\bibitem{fw2}D. Fang and C. Wang,
\textit{Almost global existence for some semilinear wave equations with almost critical regularity},
Comm. Partial Differential Equations {\bf 38} (2013), 1467--1491.

%
 \bibitem{ozawa}K. Fujiwara, V. Georgiev and T. Ozawa,
 \textit{On global well-posedness for nonlinear semirelativistic equations in some scaling subcritical and critical cases},
 arXiv:1611.09674.
%
\bibitem{guo} Z. Guo, J. Li, K. Nakanishi, L. Yan,
 \textit{On the boundary Strichartz estimates for wave and Schrodinger equations},
 arXiv:1805.01180.
 %
 \bibitem{pg}P. Grisvard,
\textit{ Elliptic problems in nonsmooth domains},
 Monographs and Studies in Mathematics, 24. Pitman (Advanced Publishing Program), Boston, MA, 1985.
 %
  \bibitem{H-W}K. Hidano and C. Wang,
 \textit{Fractional derivatives of composite functions and the Cauchy problem for the nonlinear half wave equation},
  arXiv:1707.08319.
  %
\bibitem{hjlw} K. Hidano, J. Jiang, S. Lee, C. Wang
\textit{Weighted fractional chain rule and nonlinear wave equations with minimal regularity},
arXiv:1605.06748.
 %
 \bibitem{Kl80}S. Klainerman,
\textit{Global existence for nonlinear wave equations},
Comm. Pure Appl. Math.  (1980), 43--101.
%
 
\bibitem{km}S. Klainerman and M. Machedon,
\textit{Estimates for null forms and the spaces $H_{s,\delta}$},
Internat. Math. Res. Notices  (1996), 853--865.
%
\bibitem{km1}S. Klainerman and M. Machedon,
\textit{On the algebraic properties of the $H_{n/2,1/2}$ spaces},
Internat. Math. Res. Notices (1998), 765--774.

\bibitem{ks}S. Klainerman and S. Selberg,
\textit{ Bilinear estimates and applications to nonlinear wave equations},
 Commun. Contemp. Math. {\bf 2} (2002), 223--295.
%
\bibitem{lbd3}H. Lindblad,
\textit{A sharp counterexample to the local existence of low-regularity solutions to nonlinear wave equations},
Duke Math. J. {\bf 72} (1993), 503--539.
%
\bibitem{lbd}H. Lindblad,
\textit{Counterexamples to local existence for semi-linear wave equations},
 Amer. J. Math. {\bf 118} (1996), 1--16.
 %
\bibitem{lbd2}H. Lindblad,
\textit{Counterexamples to local existence for quasi-linear wave equations},
 Math. Res. Lett. {\bf 5} (1998),  605--622.
%
\bibitem{mnno}S. Machihara, M. Nakamura, K. Nakanishi, T. Ozawa,
\textit{Endpoint Strichartz estimates and global solutions for the nonlinear Dirac equation},
J. Funct. Anal. {\bf 219} (2005), 1--20.




\bibitem{ps}G. Ponce and T. C. Sideris,
\textit{Local regularity of nonlinear wave equations in three space dimensions},
Comm. Partial Differential Equations. {\bf  18} (1993),  169--177.
%

%
\bibitem{ster}J. Sterbenz,
\textit{Global regularity and scattering for general non-linear wave equations. II. (4+1) dimensional Yang-Mills equations in the Lorentz gauge},
 Amer. J. Math. {\bf 129}, (2007), 611--664.

%
\bibitem{Tao}T. Tao,
\textit{ Nonlinear dispersive equations, local and global analysis},
CBMS Regional Conference Series in Mathematics, 106. American Mathematical Society, Providence, RI (2006).
%
\bibitem{dt}D. Tataru,
\textit{On the equation $\Box u = |\nabla u|^{2}$ in 5+1 dimensions},
Math. Res. Lett. {\bf 6} (1999), 469--485.
%
\bibitem{zhou}Y. Zhou,
\textit{On the equation $\Box u = |\nabla u|^{2}$ in four space dimensions},
 Chinese Ann. Math. Ser. {\bf 24} (2003), 293--302.


%
\end{thebibliography}

\end{document}